# Robust optimal component design under consideration of local material defects

Jannis Greifenstein, Mathematical Optimization, Collaborative Research Center 814 - Additive Manufacturing, Friedrich-Alexander-Universität Erlangen-Nürnberg, Germany, jannis.greifenstein@fau.de

Michael Stingl, Mathematical Optimization, Collaborative Research Center 814 – Additive Manufacturing, Friedrich-Alexander-Universität Erlangen-Nürnberg, Germany, michael.stingl@fau.de

*Abstract*—An important issue in additive manufacturing is the reliability and reproducibility of parts. One major problem in achieving this are uncontrolled local variations in the obtained material properties which arise in the complex manufacturing process and are usually not taken into account in the design of components.

We consider the optimal layout of a part to withstand a given loading, under the assumption that the local material properties are not precisely known. The material uncertainties are treated by a worst case approach. This means that for each layout a given amount of defects in material properties is distributed in the design domain, such that the stiffness of the component is maximally weakened. As a consequence, an optimization result is obtained which is as insensitive as possible with respect to unknown variations in the material parameters.

The general model is introduced and an algorithm for its solution using gradient based methods is suggested. Finally, numerical results are presented and discussed.

1. INTRODUCTION

In recent years, Additive Manufacturing (AM) processes have received ever-increasing attention. Compared with more traditional methods, AM allows for almost arbitrary design freedom in the layout of components and is subject to relatively little manufacturing constraints. On the other hand, the manufacturing process is very complex and still lacks in reproducibility of parts. This is not only due to manufacturing precision, but also a result of varying inhomogeneities in the material properties. These inhomogeneous material properties can be already caused by small thermal imbalances in the build chamber or simply variations in powder quality. Even running the same build twice may result in two mechanically different components.

To evaluate the mechanical performance of a component, we will consider the compliance functional, which is a good measure for a structure's stiffness with respect to a prescribed loading. The goal will be to find a structure with the smallest compliance, i.e. with the largest structural stiffness. The purpose of this work is to methodically analyze the influence that uncertainties in the local material properties can have on the performance of a component. Furthermore, we look for an optimization result, which is as insensitive to these uncertain material properties as possible. In the literature, different formulations for robust topology optimization can be found. Most publications consider a robustness with respect to uncertainties in the design layout itself, such as in [1, 2]. One of the few publications considering uncertainties in the local material properties is [3], where the uncertainties are treated as stochastic fields. In contrast, in this work we will consider the maximum influence uncertain material properties can have on the structural performance using a non-probabilistic worst-case approach, which is more related to works in optimization with damage such as [4]. This means that simultaneously to the compliance minimization, the maximal negative effect is evaluated that a certain amount of material defects can have on component performance. To find the worst-case distribution of material defects, we use a second optimization problem. Mathematically, this leads to a so-called MinMax problem. The compliance is minimized with respect to (w.r.t.) the layout of the structure and maximized w.r.t. the material defects that might occur. MinMax problems are usually non-smooth problems and very hard to solve. This means that we will have to choose a model which still allows a numerical solution in a reasonable amount of time.

## 2. PROBLEM DESCRIPTION

### 2.1. State problem

The state problem of linear elasticity is discretized using the finite element method (FEM), yielding the following linear system of equations (see e.g. [5]):

$$Ku = f, \qquad (1)$$

where u is the finite element displacement vector, f is the vector with the applied nodal forces and the finite element stiffness matrix K is assembled as

$$K = \sum_{e \in \Omega} K_e, \quad K_e = \sum_{l=1}^{n_g} B_{el}^{\mathrm{T}} C_e B_{el}. \qquad (2)$$

Here, $\Omega$ denotes the simulation domain which is partitioned into finite elements $e$, $B$ denote strain displacement matrices, $n_g$ the number of integration points and $C$ the linear elastic stiffness tensor from Hooke's law. Note that $C_e$ will be different from element to element in the optimization setting according to the material model described in the next paragraph.

### 2.2. Material model

In practice, the observed variations in material properties are larger in build direction than within the building layer. Nevertheless, we choose an isotropic model to keep the final model numerically manageable. We base the chosen formulation on a publication about truss topology design with degradation effects [4]. The formula for the linear elastic tensor in Hooke's law is given as

$$\hat{C}_e(\delta) = \left(\frac{1-\delta_e}{E_0} + \frac{\delta_e}{E_D}\right)^{-1} \frac{1}{(1+\nu)(1-2\nu)} \begin{pmatrix} 1-\nu & \nu & \nu & 0 & 0 & 0 \\ \nu & 1-\nu & \nu & 0 & 0 & 0 \\ \nu & \nu & 1-\nu & 0 & 0 & 0 \\ 0 & 0 & 0 & (1-2\nu)/2 & 0 & 0 \\ 0 & 0 & 0 & 0 & (1-2\nu)/2 & 0 \\ 0 & 0 & 0 & 0 & 0 & (1-2\nu)/2 \end{pmatrix}, \qquad (3)$$

where $e$ denotes the finite element number, $E_0$ and $E_D$ the largest and smallest Young's moduli occurring in the fluctuations and $\nu$ is the Poisson's ratio of the material. The parameter $\delta$ is used to interpolate between the least stiff and the stiffest material. This parameter represents the locally varying material properties. A major advantage of this formulation will be the concavity of the compliance functional with respect to $\delta$ in the optimization problem introduced in the next subsection. An example of the effective material stiffness for the different values of $\delta$ is shown in Fig. 1.

The distribution of the material for the actual structure is represented through a continuous so-called pseudo-density $\rho$. The pseudo-density $\rho$ attains values between a small positive number $\rho_{min}$ (representing void) and 1 (representing a material point). To avoid obtaining a lot of physically meaningless intermediate densities, the pseudo-density is penalized using the well-known SIMP approach [6]. In SIMP, the material tensor is multiplied by a penalized density $\rho^p$, effectively reducing the stiffness of elements with intermediate densities $\rho_e$ for $p$ larger than 1. Consequently, the formula for the stiffness tensor is

$$C_e(\rho, \delta) = \rho_e^p \hat{C}_e(\delta) \qquad (4)$$

with $\hat{C}_e$ as in (3).

### 2.3. Optimization problem

In the optimization problem, the compliance $f^{\mathrm{T}} u$ is minimized with respect to the design layout $\rho$ (the structural stiffness is maximized). At the same time, the uncertain material parameters $\delta$ are subject to a worst-case approach, i.e. the compliance is maximized w.r.t. $\delta$. The resulting MinMax problem is

$$\begin{cases} \min_{\rho \in U_{ad}} \max_{\delta \in \Delta_{ad}} f^{\mathrm{T}} u(\rho) \\ \text{s.t.: } K(\rho, \delta) u = f, \end{cases} \qquad (5)$$

where $U_{ad}$ and $\Delta_{ad}$ are the sets of admissible designs and $K(\rho, \delta)$ is assembled as in (2) with $C_e$ from (3). A common problem for this SIMP model is the occurrence of checkerboards. To avoid this, we use a density filter as introduced in [7, 8] for the density variable $\rho$. As admissible design sets, we choose

$$U_{ad} = \left\{ \rho \in \mathbb{R}^n : 0 < \rho_{min} \leq \rho_e \leq 1, \; e = 1, \ldots, n, \; \sum_{e \in \Omega} v_e \rho_e / |\Omega| \leq V \right\},$$

$$\Delta_{ad} = \left\{ \delta \in \mathbb{R}^n : 0 \leq \delta_e \leq 1, \; e = 1, \ldots, n, \; \sum_{e \in \Omega} v_e \rho_e^p \delta_e / |\Omega| = 0.5, \; \sum_{e \in \Omega} v_e (\delta_e - 0.5)^2 / |\Omega| = D \right\},$$

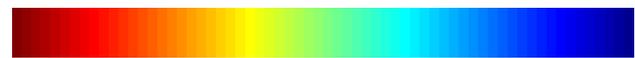

(a) $\delta$ ranging from 1 (left) to 0 (right)

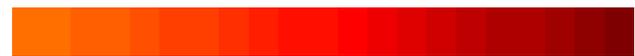

(b) Effective material coefficient for $\delta$ as shown above (left: $E_D = 0.75$, right: $E_0 = 1$)

Figure 1.   Illustration of the material model

where $v_e$ is the volume of finite element $e$, $V$ is the admissible volume fraction of material and $D$ is the allowed fraction of material defects. The choice of $\Delta_{ad}$ causes that the material defects will on average yield the expected value in the middle between the stiffest and the least stiff material. Furthermore, the constraint on the amount of material defects is summed up quadratically, such that material at the average of 0.5 is "free" and material at the margins of the allowed range is costlier.

MinMax problems like problem (5) are usually non-smooth problems. Thus, they require special solution algorithms which are generally not capable of treating the amount of optimization variables arising in structural optimization. Consequently, in the following we modify the optimization problem to allow for a more efficient numerical solution.

We start by considering an equivalent problem to (5) by using the principle of potential energy (i.e. we replace the state problem with a minimization in the displacement vector $u$):

$$\left\{ \min_{\rho \in U_{ad}} \max_{\delta \in \Delta_{ad}, u \in V} 2f^\top u - \sum_{e=1,\ldots,n} u^\top K_e(\rho,\delta) u, \right. \quad (6)$$

where $V$ is the discretized space of admissible displacements.

In the following, we will consider (6) as a bilevel optimization problem and write equivalently

$$\min_{\rho \in U_{ad}} f^\top u(\rho) \quad (7)$$

$$\text{s.t.: } (\delta, u)(\rho) = \operatorname*{argmin}_{\delta, u \in V} -2f^\top u + \sum_{e=1,\ldots,n} u^\top K_e(\rho,\delta) u,$$

$$0 \le \delta_e \le 1,\ e = 1,\ldots,n,\ \sum_{e \in \Omega} v_e \rho_e^p \delta_e / |\Omega| = 0.5,$$

$$\sum_{e \in \Omega} v_e (\delta_e - 0.5)^2 / |\Omega| = D.$$

Looking at the optimality conditions of the inner minimization problem in (7), we see that the inequality bound constraints on $\delta$ are the only reason the optimality conditions are not a system of equations. If it were a set of equations, the MinMax problem would actually be differentiable and a standard adjoint calculus with the optimality system of the inner minimization problem could be used. This lead us to apply a barrier method to the bound constraints, effectively smoothing the problem:

$$\min_{\rho \in U_{ad}} f^\top u(\rho) \quad (8)$$

$$\text{s.t.: } (\delta, u)(\rho) = \operatorname*{argmin}_{\delta, u \in V} -2f^\top u + \sum_{e=1,\ldots,n} u^\top K_e(\rho,\delta) u$$

$$- \mu^* \sum_{e=1,\ldots,n} \left( \log(\delta_e) + \log(1-\delta_e) \right),$$

$$\sum_{e \in \Omega} v_e \rho_e^p \delta_e / |\Omega| = 0.5,\ \sum_{e \in \Omega} v_e (\delta_e - 0.5)^2 / |\Omega| = D$$

where $\mu^*$ is a small positive barrier parameter.

While the cost functional of the inner problem was already convex simultaneously in $\delta$ and $u$ before, this modified cost functional is now even strictly convex. Furthermore, both formulations have a sparse Hessian and only two constraints, allowing for the efficient use of second order algorithms.

3. NUMERICAL RESULTS

As solution algorithm, we solve the inner minimization problem in $\delta$ and $u$ with the second order interior point optimizer IPOPT [9], then compute the gradients w.r.t. using an adjoint calculus and and optimize $\rho$ with the so-called Method of Moving Asymptotes (MMA, see [10]). This procedure is repeated until convergence is achieved.

3.1. Loaded cantilever

As exemplary setup we use a cantilever fixed on the left side with a line load at the right bottom: The domain is discretized using 400×200 rectangular finite elements. An illustration of the setup can be found in Fig. 2.

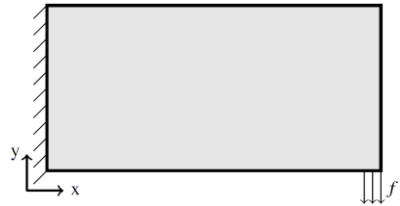

Figure 2. Visualization of the exemplary setup

We use a penalization factor of p = 5 and a linearly decaying density filter with a radius corresponding to 7 element lengths. As resource constraints, we use a total relative material volume of 50% (V = 0.5). The uncertain material properties are allowed to deviate by 25%, i.e. the normalized values in the optimization problem are $E_D = 0.75$ and $E_0 = 1$. These values are based on tensile testing results for AlSi10Mg [11]. Note that the variations in Young's modulus will change for different processes, materials, machines and process parameters. E.g. for reused polyamide powder, the elastic modulus may decrease rapidly up to half the original value after several uses [12] while for the maraging steel 18Ni-300 the variations may be close to 14% (see [13]).

As it is difficult to know the actual amount of variations occurring in the material, we consider different fractions for this constraint. Specifically, we use D = 0.01, D = 0.02 and D = 0.04, i.e. a maximum of 4%, 8% and 16% of $\delta$ could be at the upper and lower bounds.

In the results shown, the result of a standard SIMP optimization with fixed material properties is used as the starting value in the proposed algorithm. As a reference value, we use the compliance of this standard topology optimization result with the average material of $\delta = 0.5$ on the whole domain. We compare this value with the compliance of the same part with the worst-case distribution of material defects and with the worst-case compliance of the final robust optimization result. The results are visualized in Fig. 3 and 4 and the relative compliance values can be found in Table 1.

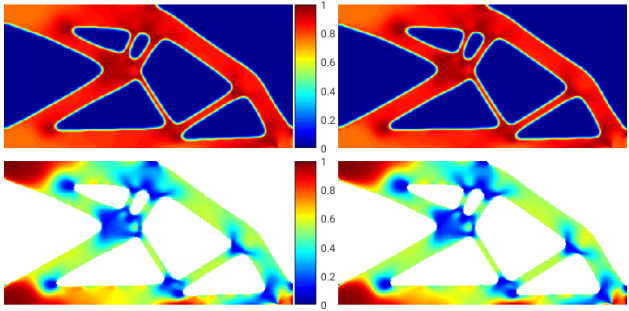

Figure 3. Visualization of the optimization results for $D = 0.02$. Upper row: physical stiffness with worst-case approach. Bottom row: worst-case distribution of material defects on topology thresholded at physical density 0.4. Left column: standard SIMP optimization result. Right column: result of robust optimization

TABLE I.   COMPLIANCE VALUES (REFERENCE AND WORST-CASE)

|  | SIMP reference ($\delta \equiv 0.5$) | SIMP | robust |
|---|---|---|---|
| $D = 0.01$ | 100% | 102.31% | 102.16% |
| $D = 0.02$ | 100% | 102.83% | 102.61% |
| $D = 0.04$ | 100% | 103.31% | 102.89% |

These first results suggest that the structures optimized with the standard SIMP approach seem to already be close to the robust optimization result. The actual influence of the material defects in this example is even with the worst-case approach smaller than expected. This by itself is already an interesting find. Furthermore, the influence is noticeable, but small for the normal topology optimization result considered here.

4. CONCLUSION

We proposed a model and an algorithm for the robust optimal design of loaded components with unpredictable local material variations and applied it successfully to a numerical example. Specifically the solution approach for the inner problem proved to be very efficient. After the first view iterations it only needed 5 iterations for solving the inner optimization problem to a KKT precision of less than $10^{12}$.

The preliminary results suggest, that at least for the compliance functional, the specific chosen example and an already topology optimized part, the influence of the uncontrolled material defects is not negligible, but still minor. The robust optimization led to better results compared with the standard SIMP approach, however not by much.

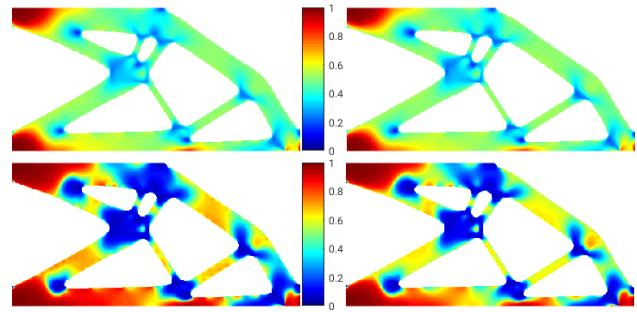

Figure 4. Visualization of the optimization results for $D = 0.01$ (upper row) and $D = 0.04$ (bottom row). For an explanation of the visualization see Fig. 3 for the bottom row

Future research could include different numerical examples to see whether the preliminary findings hold up. Furthermore, a better description of the uncertainties in the material properties observed in reality is planned. Finally, the proposed smoothed optimization variant may be suitable for other, more difficult cost functionals than the compliance.

ACKNOWLEDGMENT

The authors want to thank the German Research Foundation (DFG) for funding Collaborative Research Centre 814 (CRC 814), sub-project C2.